\documentclass[11pt, reqno]{amsart}
\usepackage{amsmath, amsthm, amscd, amsfonts, amssymb, graphicx, color}
\usepackage[bookmarksnumbered, colorlinks, plainpages]{hyperref}
\hypersetup{colorlinks=true,linkcolor=red, anchorcolor=green, citecolor=cyan, urlcolor=red, filecolor=magenta, pdftoolbar=true}

\newtheorem{theorem}{Theorem}[section]
\newtheorem{lemma}[theorem]{Lemma}
\newtheorem*{theorem*}{Question}
\newtheorem{proposition}[theorem]{Proposition}
\newtheorem{corollary}[theorem]{Corollary}
\theoremstyle{definition}
\newtheorem{definition}[theorem]{Definition}
\newtheorem{example}[theorem]{Example}

\theoremstyle{remark}
\newtheorem{remark}[theorem]{Remark}
\numberwithin{equation}{section}

\begin{document}

\setcounter{page}{1}


\title[A Family of Semi-norms in $C^*$-algebras] { A Family of Semi-norms in $C^*$-algebras }
\author[Athul Augustine, Pintu Bhunia \MakeLowercase{and} P. Shankar]{Athul Augustine, Pintu Bhunia \MakeLowercase{and} P. Shankar}

\address{Athul Augustine, Department of Mathematics, Cochin University of Science And Technology,  Ernakulam, Kerala - 682022, India. }
\email{\textcolor[rgb]{0.00,0.00,0.84}{athulaugus@gmail.com, athulaugus@cusat.ac.in}}

\address{Pintu Bhunia, Department of Mathematics, Indian Institute of Science, Bengaluru 560012,
Karnataka, India. }
\email{\textcolor[rgb]{0.00,0.00,0.84}{pintubhunia5206@gmail.com, pintubhunia@iisc.ac.in}}

\address{P. Shankar, Department of Mathematics, Cochin University of Science And Technology,  Ernakulam, Kerala - 682022, India.}
\email{\textcolor[rgb]{0.00,0.00,0.84}{shankarsupy@gmail.com, shankarsupy@cusat.ac.in}}

\subjclass[2020]{47A12, 47A30, 26E60, 46L05}

\keywords{Semi-norm,  Algebraic numerical radius, Norm inequality, Mean}


\begin{abstract}
We introduce a new family of non-negative real-valued functions on a $C^*$-algebra $\mathcal{A}$, i.e., for $0\leq \mu \leq 1,$
$$\|a\|_{\sigma_{\mu}}= \text{sup}\left\lbrace \sqrt{|f(a)|^2 \sigma_{\mu} f(a^*a)}: f\in \mathcal{A}', \, f(1)=\|f\|=1 \right\rbrace, \quad $$
where $a\in \mathcal{A}$ and $\sigma_{\mu}$ is an interpolation path of the symmetric mean $\sigma$. 
These functions are semi-norms as they satisfy the norm axioms, except for the triangle inequality. Special cases satisfying triangle inequality, and a complete equality characterization is also discussed. Various bounds and relationships will be established for this new family, with a connection to the existing literature in the algebra of all bounded linear operators on a Hilbert space.
\end{abstract}
\maketitle

\section{Introduction}
 Let $\mathcal{A}$ be a unital $C^*$-algebra with unit $1$ and $\mathcal{A}'$ denotes the space of all continuous linear functionals on $\mathcal{A}$. The  set of normalized states is defined by $\mathcal{S}(\mathcal{A}) = \{f\in \mathcal{A}': f(1)=\|f\|=1\}$. The \textit{algebraic numerical range} of an element $a\in \mathcal{A}$ is defined by
$$V(a):=\{f(a): f\in \mathcal{S}(\mathcal{A})\},$$ 
which is a compact and convex subset of the complex plane $\mathbb{C}$.
The \textit{algebraic numerical radius} of $a\in \mathcal{A}$ is defined by 
$$v(a):= \sup \{|z| : z\in V(a)\}. $$
The algebraic numerical radius defines an equivalent norm on $\mathcal{A}$ via the relation
$$\frac{1}{2}\|a\|\leq v(a)\leq \|a\| \quad \forall a\in \mathcal{A}.$$
These inequalities are sharp;
 $v(a)=\|a\|$ if $a$ is normal and $v(a)=\frac{\|a\|}{2}$  if $a^2=0$ (see \cite{bhunia2024inequalities}). 
 
  If $\mathcal{A} = B(\mathcal{H})$ is the algebra of all bounded linear operators on a complex Hilbert space $\mathcal{H}$ and $A\in B(\mathcal{H})$, then $V(A) = \overline{W(A)},$ where $W(A) := \{\langle Ax,x \rangle : x\in \mathcal{H}~\text{and} ~ \|x\|=1\}$  is the \textit{numerical range} of $A$.  Consequently, $v(A)$ coincides $w(A)$,  where  $w(A) := \sup \{|z| : z\in W(A)\}$ is the \textit{numerical radius} of $A$.

A mean is a non negative function $\sigma:[0,\infty)\times [0,\infty)\rightarrow [0,\infty)$ that satisfies the following conditions (see \cite{bhatia, geometric}):
\begin{enumerate}
\item[(i)]$\sigma(a,b)\geq 0;$
\item[(ii)]if $a\leq b$, $a\leq \sigma(a,b)\leq b;$
\item[(iii)]$\sigma(a,b)$ is monotone increasing in both a and b;
\item[(iv)]For $\alpha>0$, $\sigma(\alpha a,\alpha b)=\alpha\sigma(a,b)$ (homogeneity);
\item[(v)] $\sigma(a,b)$ is continuous.
\end{enumerate}

For simplicity, we use $a\sigma b$ instead of $\sigma(a,b)$.
$\sigma$ is called a symmetric mean, if $a\sigma b= b\sigma a$.  For a symmetric mean $\sigma$, a parameterized operator mean $\sigma_\mu$, for each $\mu\in[0,1]$, is called an interpolation path \cite[Section 5.3]{mond} for $\sigma$ if it satisfies the following conditions:
\begin{enumerate}
\item[(i)]$a\sigma_0 b=a,a\sigma_1 b=b,a \sigma_\frac{1}{2} b=a\sigma b~ \forall ~ a,b\geq0;$
\item[(ii)]$(a\sigma_\mu b)\sigma(a\sigma_\nu b) = a\sigma_\frac{\mu+\nu}{2} b, \forall~\mu,\nu \in [0,1];$
\item[(iii)]For each $0\leq \mu \leq 1$, $\sigma_\mu$ is increasing in each of its components.
\end{enumerate}

A well known example of a mean is the arithmetic mean. The arithmetic mean of two non-negative numbers $a,b$ is defined by $a\nabla b=\frac{a+b}{2}$. This is a symmetric mean and the interpolation path for $\nabla$ is defined by
$$ a\nabla_\mu b= (1-\mu)a +\mu b,\quad 0\leq \mu\leq 1.$$

Other means familiar to the reader are the geometric mean and the harmonic mean of two non-negative numbers $a,b$ are defined respectively as $a\sharp b =\sqrt{ab}$ and $a ! b= \left(\frac{a^{-1}+b^{-1}}{2}\right)^{-1}$. The interpolation paths for these symmetric means are defined by
$$a\sharp_\mu b = a^{1-\mu}b^\mu, \quad \text{where $0\leq \mu\leq 1$}$$
and
$$a !_\mu b= \left((1-\mu)a^{-1} +\mu b^{-1}\right)^{-1}, \quad \text{where $0\leq \mu\leq 1.$}$$

The study of means has been extended from non-negative numbers to positive operators on Hilbert spaces.

An interpolation path of a mean that lies between the numerical radius and the operator norm was introduced in \cite{semi}.
\begin{definition}\cite[Definition 1.1]{semi}
Let $T\in B(\mathcal{H})$ and let $\sigma_{\mu}$ be an interpolation path of the symmetric mean $\sigma$. We define
$$\|T\|_{\sigma_{\mu}}= \text{sup}\left\lbrace \sqrt{|\langle Tx,x\rangle |^2 \sigma_{\mu} \|Tx\|^2}: x\in \mathcal{H}, \|x\|=1\right\rbrace, \quad 0\leq \mu \leq 1. $$
\end{definition}

These quantities are semi-norms as it satisfies the norm axioms, except the triangle inequality. This newly introduced quantity was explored in depth, with various bounds and special properties being characterized. In particular, the quantity $\|T\|_{\nabla}$ was studied extensively and $\|\cdot\|_{\nabla}$ defines a norm on $B(\mathcal{H})$.

In this work, we extend this definition to general  $C^*$-algebra. We introduce new numerical quantities that lie between the algebraic numerical radius and the $C^*$-algebra norm.
\begin{definition}
Let $a\in \mathcal{A}$ and let $\sigma_{\mu}$ be an interpolation path of the symmetric mean $\sigma$. We define
$$\|a\|_{\sigma_{\mu}}= \text{sup}\left\lbrace \sqrt{|f(a)|^2 \sigma_{\mu} f(a^*a)}: f\in \mathcal{S}(\mathcal{A})\right\rbrace, \quad 0\leq \mu \leq 1. $$
\end{definition}

If $f$ is a positive linear functional on $\mathcal{A}$, then
$ |f(a)|^2 \leq \|f\|f(a^*a).$
Therefore, for $f \in \mathcal{S}(\mathcal{A})$, we have $|f(a)|^2 \leq f(a^*a).$
This implies  
$$ v(a) \leq \|a\|_{\sigma_{\mu}} \leq \|a\|, \quad  \forall \, 0\leq \mu \leq 1. $$

\begin{example}
Let $M_2(\mathbb{C})$ denote the set of all $2\times 2$ matrices with complex entries.
Consider the matrix $A\in M_2(\mathbb{C})$:
$$A =\begin{bmatrix} 
	0 & 2 \\
	0 & 0\\
	\end{bmatrix}.$$
	We have $\|A\|=2$ and the numerical range of $A$ is the disc with center at the origin and radius $1$. Therefore, the algebraic numerical radius of $A$, $v(A)=w(A)=1.$
	Let $\sigma$ be the arithmetic mean and let $\mu=\frac{1}{2}$. Then we have $\sigma_{\mu} =\nabla$ and
	 $$\|A\|_{\nabla} = \text{sup}\left\lbrace \sqrt{\frac{|\langle Ax,x\rangle|^2 + \|Ax\|^2}{2}}:x\in \mathbb{C}^2, \|x\|=1\right\rbrace= \sqrt{\frac{3}{2}}.$$
Therefore,  we have
$1=v(A)< \sqrt{\frac{3}{2}}= \|A\|_{\nabla} < \|A\|=2.$
\end{example}

The main aim of this work is to study this new quantity, and deduce bounds and relations in connection to the existing literature of operators.

We end this section by noting the following definitions of some special classes of elements in the $C^*$-algebra that will be used later. An element $a\in \mathcal{A}$ is positive (we write $a\geq 0$) if $a$ is self-adjoint and the spectrum of $a$ is contained in $\mathbb{R}^+$. For an arbitrary element $a\in \mathcal{A}$, $a^*a$ is always positive.
 An element $a$ is called $p$-hyponormal with $0<p\leq 1$ if $(a^*a)^p - (aa^*)^p \geq 0$. If $p=1$, a is hyponormal and if $p=\frac{1}{2}$, then a is semi-hyponormal. For real nmbers $\alpha$ and $\beta$ with $0<\alpha\leq 1\leq \beta,$ $a$ is called $(\alpha,\beta)$-normal \cite{alpha}, if $\alpha^2|a|^2\leq |a^*|^2\leq \beta^2|a|^2$, where $|a|^2=a^*a$.
 
\section{Basic properties}

In this section, we discuss some basic properties of $\|\cdot\|_\sigma$. The following proposition can be deduced directly from the definition of $\|\cdot\|_{\sigma_\mu}$.
\begin{proposition}\label{prop2.1}
Let $a\in \mathcal{A}$ and let $\sigma_{\mu}$ be an interpolation path of the symmetric mean $\sigma$.  Then the following results hold:

\begin{enumerate}
\item[(i)] $\|a\|_{\sigma_0} = v(a)$, $\|a\|_{\sigma_1} = \|a\|$.
\item[(ii)] For each $0\leq \mu \leq 1$, $ v(a) \leq \|a\|_{\sigma_{\mu}} \leq \|a\|.$
\item[(iii)] For each $0\leq \mu \leq 1$, $\|a\|_{\sigma_{\mu}} \leq \sqrt{v(a)^2\sigma_{\mu} \|a\|^2}$.
\item[(iv)] $\|a\|_{\sigma_{\mu}} = 0$, for any $\mu$, if and only if $a=0$.
\item[(v)] $\|\lambda a\|_{\sigma_{\mu}} = |\lambda| \|a\|_{\sigma_{\mu}}$ for any $\lambda \in \mathbb{C}$.
\item[(vi)] If $a$ is normal, then
$ v(a) = \|a\|_{\sigma_{\mu}} = \|a^*\|_{\sigma_{\mu}} = \|a\|.$
\end{enumerate}
\end{proposition}

For an element $a$ in a $C^*$-algebra $\mathcal{A}$, it is known that $\|a\|=\|a^*\|$, $\|a\|=\left\| |a| \right \|$, $v(a^*)=v(a)$ and $v(a)\leq v(|a|)$. In the next proposition, we characterize the relations for $\|\cdot\|_{\sigma_\mu}$ in this context.
\begin{proposition}
Let $a\in \mathcal{A}$ and let $\sigma_{\mu}$ be an interpolation path of the symmetric mean $\sigma$. Then the following results hold:
\begin{enumerate}
\item[(i)] If $a$ is semi-hyponormal, then $ \|a\|_{\sigma_{\mu}} \leq \|~|a|~\|_{\sigma_{\mu}}$.
\item[(ii)] If $a$ is hyponormal, then $ \|a^*\|_{\sigma_{\mu}} \leq \|a\|_{\sigma_{\mu}}$.
\item[(iii)] If $a$ is $(\alpha,\beta)$-normal, then $ \alpha\|a\|_{\sigma_{\mu}}\leq \|a^*\|_{\sigma_{\mu}} \leq \beta\|a\|_{\sigma_{\mu}}$.
\end{enumerate}
\end{proposition}
\begin{proof}
$(i)$ Let $\mathcal{A}$ be a $C^*$-algebra and $f$ be a positive linear functional on $\mathcal{A}$. Then by GNS construction, for any $a \in \mathcal{A}$, we have 
$f(a)= \langle \pi_f(a)\xi,\xi\rangle,$
where $\pi_f(a)$ is the representation of the element $a$ as a bounded operator on the corresponding Hilbert space $\mathcal{H}_f$. Since $a$ is semi-hyponormal, it follows that $|a^*| \leq |a|$, and hence
\begin{equation}\label{2.2}
\begin{split}
|f(a)|^2 &= |\langle \pi_f(a)\xi,\xi\rangle|^2\\
& \leq \langle |\pi_f(a)|\xi,\xi\rangle \langle |(\pi_f(a))^*|\xi,\xi\rangle\\
&= \langle \pi_f(|a|)\xi,\xi\rangle \langle \pi_f(|a^*|)\xi,\xi\rangle\\
& = f(|a|)f(|a^*|)\\
&= \left|f(|a|)\right|^2
\end{split}
\end{equation}
 
  Also, since $|a|$ is a positive operator, it is easy to observe that $|a|^* =|a|$ and 
$f(a^*a) = f(|a|^2) = f(|a\|a|) =  f(|a|^*|a|).$
Combining these equations and the monotonicity of $\sigma_\mu$, we can conclude that
$$ |f(a)|^2 \sigma_{\mu} f(a^*a) \leq \left|f(|a|)\right|^2 \sigma_{\mu} f(|a|^*|a|).$$
By taking the supremum over all $f \in \mathcal{S}(\mathcal{A})$, we get $ \|a\|_{\sigma_{\mu}} \leq \| |a| \|_{\sigma_{\mu}}$.\\
$(ii)$ Let $a\in \mathcal{A}$ be hyponormal. Then $a^*a \geq aa^*$. Since $f \in \mathcal{S}(\mathcal{A})$, we  have $f(a^*a) \geq f(aa^*)$. Combining this inequality and the monotonicity of $\sigma_\mu$, we get
$$ |f(a^*)|^2 \sigma_{\mu} f(aa^*) \leq |f(a^*)|^2 \sigma_{\mu} f(a^*a).$$
Also, since $|f(a^*)| = |f(a)|$, we can conclude that 
$ |f(a^*)|^2 \sigma_{\mu} f(aa^*) \leq |f(a)|^2 \sigma_{\mu} f(a^*a).$
Thus, if $a$ is hyponormal, then $ \|a^*\|_{\sigma_{\mu}} \leq \|a\|_{\sigma_{\mu}}$.\\
$(iii)$ Since $a$ is $(\alpha,\beta)$-normal, we have 
 $ \alpha^2(a^*a) \leq (aa^*) \leq \beta^2(a^*a). $ Considering $f$ is a positive linear functional and from the monotonicity of $\sigma_\mu$, we get
$$\alpha^2 |f(a)|^2 \sigma_{\mu} \alpha^2f(a^*a) \leq |f(a^*)|^2 \sigma_{\mu} f(aa^*) \leq \beta^2 |f(a)|^2 \sigma_{\mu} \beta^2f(a^*a).$$
Since $|f(a^*)| = |f(a)|$, we can conclude that
$ \alpha\|a\|_{\sigma_{\mu}}\leq \|a^*\|_{\sigma_{\mu}} \leq \beta\|a\|_{\sigma_{\mu}}.$
\end{proof}

One can easily deduce from Proposition \ref{prop2.1} that $\|\cdot\|_{\sigma_\mu}$ is a semi-norm as it satisfies the norm properties, except the triangle inequality. Here we consider a special case, $\sigma=\nabla$ to get the following result.
\begin{proposition}\label{prop2.4}
Let $a,b \in \mathcal{A}$. Then 
$$\|a+b\|_{\nabla} \leq \|a\|_{\nabla} + \|b\|_{\nabla}.$$
\end{proposition}
\begin{proof}
Let $a,b \in \mathcal{A}$. Then 
$$\|a+b\|^2_{\nabla} = \text{sup}\left\lbrace \frac{1}{2}\left(|f(a+b)|^2 + f\left((a+b)^*(a+b)\right)\right): f\in \mathcal{S}(\mathcal{A})\right\rbrace.$$
For any $ f\in \mathcal{S}(\mathcal{A})$, we have
\begin{eqnarray}\label{2.2}
&& \frac{1}{2}  \left(|f(a+b)|^2 + f\left((a+b)^*(a+b)\right)\right) \notag \\
&&\leq \frac{1}{2}\left(\left(|f(a)| + |f(b)| \right)^2 + f\left((a+b)^*(a+b)\right)\right).
\end{eqnarray}
Now, since $f(a^*b) = \overline{f(b^*a)}$, we get
\begin{eqnarray}\label{2.3}
f\left((a+b)^*(a+b)\right) &=& f(a^*a) + f(a^*b) + f(b^*a) + f(b^*b) \notag \\
&=& f(a^*a) + f(b^*b) + 2\mathcal{R}(f(a^*b)),
\end{eqnarray}
where $\mathcal{R}$ denotes the real part. From the Cauchy-Schwarz inequality of positive linear functionals on a $C^*$-algebra, we have
$|f(a^*b)|^2 \leq f(a^*a)f(b^*b).$
Taking the square root on both sides, we get
$|f(a^*b)| \leq \sqrt{f(a^*a)f(b^*b)}.$
This implies that
$$2\mathcal{R}(f(a^*b)) \leq 2|f(a^*b)| \leq 2\sqrt{f(a^*a)f(b^*b)}.$$
Substituting this in \eqref{2.3}, we get
\begin{equation*}
\begin{split}
f\left((a+b)^*(a+b)\right) &\leq  f(a^*a) + f(b^*b) + 2\sqrt{f(a^*a)f(b^*b)}\\
&= \left(\sqrt{f(a^*a)} + \sqrt{f(b^*b)}\right)^2.
\end{split}
\end{equation*}
Substituting this in \eqref{2.2} and employing Cauchy-Schwarz inequality, we get
\begin{equation*}
\begin{split}
\frac{1}{2} & \left(|f(a+b)|^2 + f\left((a+b)^*(a+b)\right)\right) \\
&\leq \frac{1}{2}\left(\left(|f(a)| + |f(b)| \right)^2 + f(a^*a) + f(b^*b) + 2\sqrt{f(a^*a)f(b^*b)}\right).\\
&=\frac{1}{2}\left(|f(a)|^2 + |f(b)|^2 + f(a^*a) + f(b^*b) + 2\left(|f(a)\|f(b)| + \sqrt{f(a^*a)f(b^*b)}\right)\right).\\
&\leq \frac{|f(a)|^2 + f(a^*a)}{2} + \frac{|f(b)|^2 +f(b^*b)}{2} + 2\sqrt{\frac{|f(a)|^2 + f(a^*a)}{2}}\sqrt{\frac{|f(b)|^2 +f(b^*b)}{2}}\\
&= \|a\|^2_{\nabla} + \|b\|^2_{\nabla} +2\|a\|_{\nabla} \|b\|_{\nabla}.
\end{split}
\end{equation*}
This completes the proof.
\end{proof}

\begin{corollary}
Let $\nabla$ be the arithmetic mean. Then $\|\cdot\|_{\nabla}$ is a norm on $\mathcal{A}$.
\end{corollary}
\begin{proof}
This follows from Propositions \ref{prop2.1} and \ref{prop2.4}.
\end{proof}
Note that the norm $\|\cdot\|_{\nabla}$ does not make  $\mathcal{A}$ a $C^*$-algebra. This is due to the fact that the norm that makes $\mathcal{A}$ a $C^*$-algebra is unique, and here  $\|a\|_{\nabla} \leq  \|a\|$.

In the next theorem, we give a characterization for the equality  $\|a+b\|_{\nabla} = \|a\|_{\nabla} + \|b\|_{\nabla}$ to hold in $\mathcal{A}$.
\begin{theorem}
Let $a,b \in \mathcal{A}$. Then the following conditions are equivalent:
\begin{itemize}
\item[(1)] $\|a+b\|_{\nabla} = \|a\|_{\nabla} + \|b\|_{\nabla}.$
\item[(2)] There exists a sequence $\{f_n\} \in \mathcal{S}(\mathcal{A})$ such that
$$\lim_{n \longrightarrow \infty} \mathcal{R}(f_n(b^*a) + \overline{f_n(a)}f_n(b)) = 2 \|a\|_{\nabla}\|b\|_{\nabla}.$$
\end{itemize}
\end{theorem}
\begin{proof}
Suppose that  $\|a+b\|_{\nabla} = \|a\|_{\nabla} + \|b\|_{\nabla}.$ By the hypothesis, there exists a sequence $\{f_n\} \in \mathcal{S}(\mathcal{A})$ such that
$$\lim_{n \longrightarrow \infty}\frac{1}{2}\left(|f_n(a+b)|^2 + f_n\left((a+b)^*(a+b)\right)\right) = (\|a\|_{\nabla} + \|b\|_{\nabla})^2. $$
For every $n\in \mathbb{N}$, we have
\begin{equation*}
\begin{split}
\frac{1}{2}&\left(|f_n(a+b)|^2 + f_n\left((a+b)^*(a+b)\right)\right)\\
&=\frac{1}{2}\left(|f_n(a)+f_n(b)|^2 + f_n(a^*a) + f_n(a^*b) +f_n(b^*a) +f_n(b^*b)  \right)\\
&= \frac{1}{2}\left(|f_n(a)|^2+|f_n(b)|^2 + 2\mathcal{R}(f_n(b)\overline{f_n(a)}) +  f_n(a^*a)  +f_n(b^*b) + 2\mathcal{R}(f_n(b^*a)) \right)\\
&\leq \|a\|^2_{\nabla} + \|b\|^2_{\nabla} + \mathcal{R}(f_n(b^*a) + \overline{f_n(a)}f_n(b) )\\
&\leq \|a\|^2_{\nabla} + \|b\|^2_{\nabla} + |(f_n(b^*a) + \overline{f_n(a)}f_n(b) )|\\
&\leq \|a\|^2_{\nabla} + \|b\|^2_{\nabla} + (\sqrt{f_n(b^*b)}\sqrt{f_n(a^*a)} + |\overline{f_n(a)}f_n(b)|)\\
&\leq \|a\|^2_{\nabla} + \|b\|^2_{\nabla} + 2 \sqrt{\frac{f_n(a^*a)+|f_n(a)|^2}{2}}  \sqrt{\frac{f_n(b^*b)+|f_n(b)|^2}{2}}\\
&\leq (\|a\|_{\nabla} + \|b\|_{\nabla})^2.
\end{split}
\end{equation*}
Then, taking limit $n\rightarrow \infty$, we conclude that 
$$\lim_{n \longrightarrow \infty} \mathcal{R}(f_n(b^*a) + \overline{f_n(a)}f_n(b)) = 2 \|a\|_{\nabla}\|b\|_{\nabla}.$$
Now suppose that $(2)$ holds, i.e.,  there exists a sequence $\{f_n\} \in \mathcal{S}(\mathcal{A})$ such that
$$\lim_{n \longrightarrow \infty} \mathcal{R}(f_n(b^*a) + \overline{f_n(a)}f_n(b)) = 2 \|a\|_{\nabla}\|b\|_{\nabla}.$$
Then
\begin{equation*}
\begin{split}
\mathcal{R}^2 & \left( \frac{1}{2}(f_n(b^*a) + \overline{f_n(a)}f_n(b))\right)\\
&=\left| \frac{1}{2}(f_n(b^*a) + \overline{f_n(a)}f_n(b))\right|^2 - \mathcal{I}^2\left( \frac{1}{2}(f_n(b^*a) + \overline{f_n(a)}f_n(b))\right)\\
& \leq \left| \frac{1}{2}(f_n(b^*a) + \overline{f_n(a)}f_n(b))\right|^2 \\
& \leq \left( \frac{1}{2}|f_n(b^*a)| + \frac{1}{2} |\overline{f_n(a)}\|f_n(b)|\right)^2 \\
& \leq \left( \frac{1}{2}f_n(a^*a)^{\frac{1}{2}}f_n(b^*b)^{\frac{1}{2}} + \frac{1}{2} |\overline{f_n(a)}\|f_n(b)|\right)^2 \\
& \leq \left(\frac{f_n(a^*a)+|f_n(a)|^2}{2} \times  \frac{f_n(b^*b)+|f_n(b)|^2}{2}\right)\\
& \leq \|a\|^2_{\nabla}\|b\|^2_{\nabla}.
\end{split}
\end{equation*}
Therefore, we get
\begin{equation}\label{2.4}
\mathcal{R}^2  \left( \frac{1}{2}(f_n(b^*a) + \overline{f_n(a)}f_n(b))\right)\leq \|a\|^2_{\nabla}\|b\|^2_{\nabla}.
\end{equation}
Since 
$$\lim_{n \longrightarrow \infty} \mathcal{R}(f_n(b^*a) + \overline{f_n(a)}f_n(b)) = 2 \|a\|_{\nabla}\|b\|_{\nabla}$$
and
$$\left(\frac{f_n(a^*a)+|f_n(a)|^2}{2}\right) \leq  \|a\|^2_{\nabla}, \left(\frac{f_n(b^*b)+|f_n(b)|^2}{2}\right) \leq  \|b\|^2_{\nabla},$$
we infer that 
$$\lim_{n \longrightarrow \infty}\left(\frac{f_n(a^*a)+|f_n(a)|^2}{2}\right) =  \|a\|^2_{\nabla} $$
and
$$\lim_{n \longrightarrow \infty} \left(\frac{f_n(b^*b)+|f_n(b)|^2}{2}\right) =  \|b\|^2_{\nabla}.$$
Finally, we have 
\begin{equation*}
\begin{split}
(\|a\|_{\nabla} + \|b\|_{\nabla})^2 &= \|a\|^2_{\nabla} + \|b\|^2_{\nabla} + 2\|a\|_{\nabla} \|b\|_{\nabla}\\
&=\lim_{n \longrightarrow \infty}\left(\frac{f_n(a^*a)+|f_n(a)|^2}{2}\right) + \lim_{n \longrightarrow \infty} \left(\frac{f_n(b^*b)+|f_n(b)|^2}{2}\right) \\
& \qquad +\lim_{n \longrightarrow \infty} \mathcal{R}(f_n(b^*a) + \overline{f_n(a)}f_n(b))\\
&=\lim_{n \longrightarrow \infty}\frac{1}{2}\left(|f_n(a+b)|^2 + f_n\left((a+b)^*(a+b)\right)\right)\\
& \leq (\|a + b\|_{\nabla})^2\\
& \leq (\|a\|_{\nabla} + \|b\|_{\nabla})^2,
\end{split}
\end{equation*}
where we have used Proposition \ref{prop2.4} to obtain the last inequality. Thus, we have $\|a+b\|_{\nabla} = \|a\|_{\nabla} + \|b\|_{\nabla},$ which completes the proof.
\end{proof}

\begin{remark}
The following conditions are equivalent, for $a,b \in \mathcal{A}$ and a sequence $\{f_n\} \in \mathcal{S}(\mathcal{A})$:
\begin{itemize}
\item[(1)]$\lim_{n \longrightarrow \infty} \mathcal{R}(f_n(b^*a) + \overline{f_n(a)}f_n(b)) = 2 \|a\|_{\nabla}\|b\|_{\nabla}.$
\item[(2)]$\lim_{n \longrightarrow \infty} (f_n(b^*a) + \overline{f_n(a)}f_n(b)) = 2 \|a\|_{\nabla}\|b\|_{\nabla}.$
\end{itemize}
\end{remark}
Indeed, if $(2)$ holds, then
\begin{equation*}
\begin{split}
2 & \|a\|_{\nabla}\|b\|_{\nabla}\\
&= \frac{2 \|a\|_{\nabla}\|b\|_{\nabla} + 2 \|a\|_{\nabla}\|b\|_{\nabla}}{2}\\
&= \frac{\lim_{n \rightarrow \infty} (f_n(b^*a) + \overline{f_n(a)}f_n(b)) + \lim_{n \rightarrow \infty} \overline{(f_n(b^*a) + \overline{f_n(a)}f_n(b))}}{2}\\
&= \lim_{n \longrightarrow \infty} \mathcal{R}(f_n(b^*a) + \overline{f_n(a)}f_n(b)).
\end{split}
\end{equation*}
On the other hand, assume $(1)$ holds. Then by \eqref{2.4}, we get 
$$\lim_{n \longrightarrow \infty} \mathcal{I}(f_n(b^*a) + \overline{f_n(a)}f_n(b)) = 0,$$
where $\mathcal{I}$ denotes the imaginary part and
\begin{equation*}
\begin{split}
\lim_{n \longrightarrow \infty} (f_n(b^*a) + \overline{f_n(a)}f_n(b)) &= \lim_{n \longrightarrow \infty} \mathcal{R}(f_n(b^*a) + \overline{f_n(a)}f_n(b))
= 2 \|a\|_{\nabla}\|b\|_{\nabla}.
\end{split}
\end{equation*}
\section{Some estimates on $\|a\|_{\sigma}$}
Determining the upper and lower bounds for the algebraic numerical radius and the $C^*$-algebra norm has garnered significant attention from numerous researchers. Our investigation in this section builds upon and extends these existing studies. In this section, we study lower and upper bounds for $\|a\|_\sigma$. To deduce necessary lemmas we need the following known results for Hilbert space operators.

\begin{lemma}\cite[Theorem 5]{notes}\label{2}
Let $T,S \in \mathcal{B}(\mathcal{H})$ be such that $|T|S=S^*|T|$, and let $\phi$ and $\psi$ be two non-negative continuous functions defined on $[0,\infty)$ such that $\phi(t)\psi(t)=t$ for every $t\geq 0$. Then
$$|\langle TS x,y\rangle|^2 \leq r(S)\langle \phi^2(|T|)x,x\rangle \langle \psi^2(|T^*|)y,y\rangle,$$
for every $x,y\in \mathcal{H}$, where $r(S)$ denotes the spectral radius of $S$.
In particular, if $S=I,$ the identity operator in $\mathcal{B}(\mathcal{H})$, we obtain 
$$|\langle Tx,y\rangle|^2
\leq \langle \phi^2(|T|)x,x\rangle \langle \psi^2(|T^*|)y,y\rangle,$$
for every $x,y\in \mathcal{H}$.
\end{lemma}
\begin{lemma}\cite[Theorem 1.2]{mond}\label{3}
Let $T\in \mathcal{B}(\mathcal{H})$ with spectrum in an interval $J\subseteq \mathbb{R}$. If $\phi: J\rightarrow \mathbb{R}$ is convex, then
$$\phi(\langle T x,x\rangle) \leq \langle \phi(T) x,x\rangle,$$
for all unit vector $x\in \mathcal{H}$. If $\phi$ is concave, the inequality is reversed.
\end{lemma}
Using these lemmas and GNS construction in $C^*$-algebra, we prove the following lemmas which are essential for deriving the bounds for $\|\cdot\|_\sigma.$ 

\begin{lemma}\label{3.2}
Let $a,b \in \mathcal{A}$ be such that $|a|b=b^*|a|$, and let $\phi$ and $\psi$ be two non-negative continuous functions defined on $[0,\infty)$ such that $\phi(t)\psi(t)=t$ for every $t\geq 0$. Then
$$|f(ab)|^2 \leq r(b)f(\phi^2(|a|))f(\psi^2(|a^*|)),$$
for every $f \in \mathcal{S}(\mathcal{A})$, where $r(b)$ denotes the spectral radius of $b$.
In particular, if $b=1,$ the unit element in $\mathcal{A}$, we obtain 
\begin{equation}\label{coro}
|f(a)|^2 \leq f(\phi^2(|a|))f(\psi^2(|a^*|)),
\end{equation}
for every $f \in \mathcal{S}(\mathcal{A})$.
\end{lemma}
\begin{proof}
By GNS construction, for any $a \in \mathcal{A}$, we have 
$f(a)= \langle \pi_f(a)\xi,\xi\rangle,$
where $\pi_f(a)$ is the representation of the element $a$ as a bounded operator on the corresponding Hilbert space $\mathcal{H}_f$. Then from Lemma \ref{2}, we obtain
\begin{equation*}
\begin{split}
|f(ab)|^2 &= |\langle \pi_f(ab)\xi,\xi\rangle|^2\\
& \leq r( \pi_f(b))\langle \phi^2|\pi_f(a)|\xi,\xi\rangle \langle \psi^2|(\pi_f(a))^*|\xi,\xi\rangle\\
&= r( \pi_f(b))\langle \pi_f(\phi^2(|a|))\xi,\xi\rangle \langle \pi_f(\psi^2(|a^*|))\xi,\xi\rangle\\
& \leq r(b)f(\phi^2(|a|))f(\psi^2(|a^*|)).
\end{split}
\end{equation*}
\end{proof}

\begin{lemma}\label{3.3}
Let $a\in \mathcal{A}$ with spectrum in an interval $J\subseteq \mathbb{R}$. If $\phi: J\rightarrow \mathbb{R}$ is convex, then
$\phi(f(a)) \leq f(\phi(a)), \text{ for every $f \in \mathcal{S}(\mathcal{A})$.}$
 If $\phi$ is concave, the inequality is reversed.
\end{lemma}
\begin{proof}
By GNS construction, for any $a \in \mathcal{A}$, we have 
$f(a)= \langle \pi_f(a)\xi,\xi\rangle,$
where $\pi_f(a)$ is the representation of the element $a$ as a bounded operator on the corresponding Hilbert space $\mathcal{H}_f$. Then from Lemma \ref{3}, we obtain
$$\phi(f(a)) = \phi(\langle \pi_f(a)\xi,\xi\rangle)
 \leq \langle \phi(\pi_f(a))\xi,\xi\rangle 
= \langle \pi_f(\phi(a))\xi,\xi\rangle 
 = f(\phi(a)).$$
\end{proof}

Now we are in a position to deduce the upper bound for $\|\cdot\|_\sigma.$ We derive an upper bound for the product of two elements satisfying specific commutativity conditions. 
\begin{theorem}
Let $a,b \in\mathcal{A}$ be such that $|a|b=b^*|a|$, and let $\phi$ and $\psi$ be two non-negative continuous functions defined on $[0,\infty)$ such that $\phi(t)\psi(t)=t$ for every $t\geq 0$. If $\sigma \leq \nabla$, then
$$\|ab\|_{\sigma} \leq \sqrt{\left\|\frac{r(b)}{4}(\phi^4(|a|)+\psi^4(|a^*|)) + \frac{1}{2}|ab|^2\right\|}$$
and
$$\|ab\|_{\sigma}^2 \leq \frac{1}{2}\sqrt{\left\|r(b)\phi^4(|a|)+\psi^2(|a|^2) \right\|\left\|\phi^2(|a|^2)+r(b)\psi^4(|a^*|) \right\|}.$$
\end{theorem}
\begin{proof}
Let $f\in \mathcal{S}(\mathcal{A})$. Utilizing Lemma \ref{3.2}, we can write
\begin{equation*}
\begin{split}
\sqrt{|f(ab)|^2 ~ \sigma~ f((ab)^*(ab))}
 &= \sqrt{|f(ab)|^2 ~ \sigma~ f(|ab|^2)}\\
 &\leq \sqrt{r(b)f(\phi^2(|a|))f(\psi^2(|a^*|)) ~ \sigma~ f(|ab|^2)}\\
 &\leq \sqrt{\frac{r(b)}{2} (f(\phi^2(|a|))^2+f(\psi^2(|a^*|))^2) ~ \sigma~ f(|ab|^2)}\\
 &\leq \sqrt{\frac{r(b)}{2} (f(\phi^4(|a|)+\psi^4(|a^*|))) ~ \sigma~ f(|ab|^2)}.
\end{split}
\end{equation*}
Thus we have shown that
\begin{equation}
\sqrt{|f(ab)|^2 ~ \sigma~ f((ab)^*(ab))} \leq \sqrt{\frac{r(b)}{2} (f(\phi^4(|a|)+\psi^4(|a^*|))) ~ \sigma~ f(|ab|^2)}.
\end{equation}
Using the assumption $\sigma \leq \nabla$, we get
\begin{equation*}
\begin{split}
\sqrt{|f(ab)|^2 ~ \sigma~ f((ab)^*(ab))} &\leq \sqrt{\frac{r(b)}{2} f(\phi^4(|a|)+\psi^4(|a^*|)) ~ \sigma~ f(|ab|^2)}\\
&\leq \sqrt{\frac{1}{2}\left(\frac{r(b)}{2}f (\phi^4(|a|)+\psi^4(|a^*|)) + f(|ab|^2)\right)}\\
&= \sqrt{\left(f \left(\frac{r(b)}{4}(\phi^4(|a|)+\psi^4(|a^*|)) + \frac{1}{2} |ab|^2\right)\right)}\\
&\leq \sqrt{\left\|\frac{r(b)}{4}(\phi^4(|a|)+\psi^4(|a^*|)) + \frac{1}{2} |ab|^2\right\|}.\\
\end{split}
\end{equation*}
Now, taking the supremum over all $f\in \mathcal{S}(\mathcal{A})$, we get
$$\|ab\|_{\sigma} \leq \sqrt{\left\|\frac{r(b)}{4}(\phi^4(|a|)+\psi^4(|a^*|)) + \frac{1}{2}|ab|^2\right\|}.$$
To prove the second inequality, we have

$|f(ab)|^2 ~ \sigma~ f((ab)^*(ab))$
\begin{equation*}
\begin{split}
&= |f(ab)|^2 ~ \sigma~ f(|ab|^2)\\
 &\leq r(b)f(\phi^2(|a|))f(\psi^2(|a^*|)) ~ \sigma~ \sqrt{ f(\phi ^2|ab|^2)f(\psi ^2|(ab)^*|^2)}\\
 &= \sqrt{r^2(b)(f(\phi^2(|a|)))^2)f(\psi^2(|a^*|)))^2} ~ \sigma~ \sqrt{ f(\phi ^2|ab|^2)f(\psi ^2|(ab)^*|^2)}\\
 &\leq \sqrt{r^2(b)f(\phi^4(|a|))f(\psi^4(|a^*|))} ~ \sigma~ \sqrt{ f(\phi ^2|a|^2)f(\psi ^2|a|^2)}\\
  &\leq \frac{1}{2}\left(\sqrt{r(b)f(\phi^4(|a|))r(b)f(\psi^4(|a^*|))} ~ +~ \sqrt{ f(\phi ^2|a|^2)f(\psi ^2|a|^2)}\right)\\
   &\leq \frac{1}{2}\left(\sqrt{r(b)f(\phi^4(|a|))+f(\psi ^2(|a|^2))} \sqrt{ f(\phi ^2|a|^2)+r(b)f(\psi^4(|a^*|))}\right)\\
   &\leq \frac{1}{2}\left(\sqrt{f(r(b)\phi^4(|a|)+\psi ^2(|a|^2))} \sqrt{ f(\phi ^2(|a|^2)+r(b)\psi^4(|a^*|)}\right)\\
   &\leq \frac{1}{2}\sqrt{\|r(b)\phi^4(|a|)+\psi ^2(|a|^2)\|  \|\phi ^2(|a|^2)+r(b)\psi^4(|a^*|)\|}.
\end{split}
\end{equation*}
Therefore,
$$\|ab\|_{\sigma}^2 \leq \frac{1}{2}\sqrt{\left\|r(b)\phi^4(|a|)+\psi^2(|a|^2) \right\|\left\|(\phi^2(|a|^2)+r(b)\psi^4(|a^*|)) \right\|},$$
as desired.
\end{proof}

In particular, taking $b=1$, we deduce the following upper bounds for $\|a\|_{\sigma}.$
\begin{corollary}\label{3.4}
Let $a\in \mathcal{A}$, and let $\phi$ and $\psi$ be two non-negative continuous functions defined on $[0,\infty)$ such that $\phi(t)\psi(t)=t$ for every $t\geq 0$. If $\sigma \leq \nabla,$ then
$$\|a\|_{\sigma} \leq \sqrt{\left\| \frac{1}{4}\left(\phi^4(|a|) + \psi^4(|a^*|)\right) + \frac{1}{2}|a|^2\right\|}$$
and
$$\|a\|_{\sigma}^2 \leq \frac{1}{2}\sqrt{\left\| (\phi^4(|a|) + \psi^2(|a|^2)\right\| \left\|\phi^2(|a|^2) + \psi^4(|a^*|) \right\|}.$$
\end{corollary}

In particular, for $0\leq \nu \leq 1$, if we consider $\phi(t) = t^{\nu}$ and $\psi(t) = t^{1-\nu}$ in Corollary \ref{3.4}, we obtain the following inequalities.

\begin{corollary}
Let $a\in \mathcal{A}$. Then for any $0\leq \nu \leq 1$,
$$\|a\|_{\sigma} \leq \sqrt{\left\| \frac{1}{4}\left( |a|^{4\nu} + |a^*|^{4(1-\nu)} \right) + \frac{1}{2}|a|^2\right\|}$$
and
$$\|a\|_{\sigma}^2 \leq \frac{1}{2}\sqrt{\left\| |a|^{4\nu} + |a|^{4(1-\nu)}\right\| \left\||a|^{4\nu} + |a^*|^{4(1-\nu)} \right\|}.$$
\end{corollary}

Next upper bound of $\|a\|_{\sigma}$ is as follows:

\begin{theorem}\label{alpha}
Let $a\in \mathcal{A}$. If $\sigma \leq \nabla,$ then
$$\|a\|_{\sigma}^2 \leq \frac{1}{2}\left\| (1+ \alpha) |a|^2 + (1 - \alpha) |a^*|^2\right\|,$$
for every $0\leq \alpha \leq 1$.
\end{theorem}
\begin{proof}
Let $f \in \mathcal{S}(\mathcal{A})$ and $\alpha \in [0,1]$. Then
\begin{equation*}
\begin{split}
|f(a)|^2~ \sigma~ f(a^*a) &= |f(a)|^2~ \sigma~ f(|a|^2)\\
&\leq f(|a|^{2\alpha})f(|a^*|^{2(1-\alpha)})~ \sigma~ f(|a|^2)\\
& \leq (f(|a|^{2}))^\alpha(f(|a^*|^{2})^{(1-\alpha)}~ \sigma~ f(|a|^2)\\
&\leq f(\alpha|a|^2+(1-\alpha)|a^*|^2)~ \sigma~ f(|a|^2)\\
&\leq \frac{1}{2} f((1+\alpha)|a|^2+(1-\alpha)|a^*|^2)\\
&\leq \frac{1}{2} \left\|(1+\alpha)|a|^2+(1-\alpha)|a^*|^2\right\|,
\end{split}
\end{equation*}
where we used \eqref{coro} to obtain the first inequality, Lemma \ref{3.3} to obtain the second inequality, the arithmetic-geometric mean inequality to obtain the third inequality and the fact $\sigma \leq \nabla$ to obtain the fourth inequality. Now by taking the supremum, we obtain
$\|a\|_{\sigma}^2 \leq \frac{1}{2}\left\| (1+ \alpha) |a|^2 + (1 - \alpha) |a^*|^2\right\|,$
for every $0\leq \alpha \leq 1$. 
\end{proof}

\begin{remark}
Note that for $\alpha = 0$ in Theorem \ref{alpha},  we get
$\|a\|_{\sigma}^2 \leq \frac{1}{2}\left\|  |a|^2 +  |a^*|^2\right\|.$
Therefore, the bound $v(a)\leq \|a\|_{\sigma}$ is an improvement of the existing bound \cite[Corollary 2.8]{ali}, namely,
$v^2(a) \leq \frac{1}{2}\left\|  |a|^2 +  |a^*|^2\right\|.$    
\end{remark}

The Crawford number $m(a)$, for $a\in\mathcal{A}$, is defined by the distance from the origin to the algebraic numerical range of $a$. That is $m(a) = \inf_{f\in \mathcal{S}(\mathcal{A})} |f(a)|$. In the following theorem, we obtain a lower bound for $\|a\|_{\sigma}$ in terms of $m(a)$.
\begin{theorem}\label{craw}
Let $a\in\mathcal{A}$. Then
$$\max\left( \sqrt{v^2(a)~ \sigma~m(a^*a)},\sqrt{m^2(a)~\sigma~\|a\|^2}\right) \leq \|a\|_{\sigma}.$$
\end{theorem}
\begin{proof}
Let $f\in \mathcal{S}(\mathcal{A})$, then we have
$\sqrt{|f(a)|^2~ \sigma~ f(a^*a)} \geq \sqrt{|f(a)|^2~ \sigma~ m(a^*a)}.$
Taking the supremum over all $f\in \mathcal{S}(\mathcal{A})$, we get
\begin{equation}\label{eq3.1}
\|a\|_{\sigma} \geq \sqrt{v^2(a)~ \sigma~ m(a^*a)}.
\end{equation}

On the other hand,
$\sqrt{|f(a)|^2~ \sigma~ f(a^*a)} \geq \sqrt{m^2(a)~ \sigma~ f(a^*a)}.$
Taking the supremum over all $f\in \mathcal{S}(\mathcal{A})$, we obtain
\begin{equation}\label{eq3.2}
\|a\|_{\sigma} \geq \sqrt{m^2(a)~ \sigma~ \|a\|^2}.
\end{equation}
Combining (\ref{eq3.1}) and (\ref{eq3.2}), we prove the desired result.
\end{proof}

\begin{remark}
For the case $\sigma = \nabla$, we have
$$\frac{1}{\sqrt{2}}\|a\| \leq \|a\|_{\nabla}.$$
Indeed, following Theorem \ref{craw}, we can write
\begin{equation*}
\begin{split}
\frac{1}{\sqrt{2}}\|a\| &= \max\left(\frac{1}{\sqrt{2}}v(a),\frac{1}{\sqrt{2}}\|a\|\right)\\
&\leq \max\left( \sqrt{v^2(a)~ \nabla ~m(a^*a)},\sqrt{m^2(a)~\nabla~\|a\|^2}\right)\\
& \leq \|a\|_{\nabla}.
\end{split}
\end{equation*}
\end{remark}

When $\sigma = \nabla$, Theorem \ref{craw} can be simplified as follows.
\begin{proposition}
Let $a\in\mathcal{A}$. Then
$$\max\left( \sqrt{v^2(a)m(a^*a)},\sqrt{m^2(a)\|a\|^2}\right) \leq \|a\|_{\nabla}.$$
\end{proposition}
\begin{proof}
Let $f\in \mathcal{S}(\mathcal{A})$. Then we have
\begin{equation}\label{eq3.3}
\begin{split}
\|a\|_{\nabla}^2 &\geq \frac{|f(a)|^2+ f(a^*a)}{2}\\
&= \frac{|f(a)|^2+ f(|a|^2)}{2}\\
 &\geq |f(a)| f(|a|^2)\\
 &\geq |f(a)|~m(a^*a).
\end{split}
\end{equation}
Therefore, taking the supremum over all $f\in \mathcal{S}(\mathcal{A})$, we get
$$\|a\|_{\nabla}^2  \geq v(a)~m(a^*a).$$
Also, we have
\begin{equation}\label{eq3.4}
\begin{split}
\|a\|_{\nabla}^2 &\geq |f(a)| f(|a|^2)
 \geq m(a)~ f(|a|^2).
\end{split}
\end{equation}
Taking the supremum over all $f\in \mathcal{S}(\mathcal{A})$, we get
$$\|a\|_{\nabla}^2  \geq m(a)~\|a^*a\| = m(a)~\|a\|^2.$$
Therefore, we get the desired result by combining (\ref{eq3.3}) and (\ref{eq3.4}).
\end{proof}

Finally, we provide a sufficient condition for the equality $\|  |a| |a^*| \|_{\sigma}=0$.
\begin{proposition}
Let $a\in\mathcal{A}$ be such that $a^2=0$. Then $\|  |a| |a^*| \|_{\sigma}=0$.
\end{proposition}
\begin{proof}
Let $f\in \mathcal{S}(\mathcal{A})$. Then we have
$$0\leq \sqrt{(f(|a||a^*|))^2~\sigma~ (f((|a||a^*|)^*(|a||a^*|))} \leq \sqrt{v^2(|a||a^*|)~\sigma~ \| |a||a^*| \|^2}.$$
Since $v(|a||a^*|)\leq \|~|a||a^*|~\| =\|a^2\|$, we have 
$$0\leq \sqrt{(f(|a||a^*|))^2~\sigma~ (f((|a||a^*|)^*(|a||a^*|))} \leq \sqrt{v^2(|a||a^*|)~\sigma~ \|a^2\|^2}\leq \|a^2\|.$$
Taking the supremum over all $f\in \mathcal{S}(\mathcal{A})$, we obtain
 $\|  |a| |a^*| \|_{\sigma}\leq \|a^2\|.$
 The fact $a^2=0$ will give the desired result.
\end{proof}

\textbf{Declaration of competing interest.}
There is no competing interest.

\textbf{Data availability.}
No data was used for the research described in the article.\\

{\bf Acknowledgments.} Athul Augustine is supported by the Senior Research Fellowship (09/0239(13298)/2022-EMR-I) of CSIR (Council of Scientific and Industrial Research, India). Pintu Bhunia was supported by National Post Doctoral Fellowship PDF/2022/000325 from SERB (Govt.\ of India) and SwarnaJayanti Fellowship SB/SJF/2019-20/14 from SERB (Govt.\ of India). P Shankar is supported by the Teachers Association for Research Excellence (TAR/2022/000063) of SERB (Science and Engineering Research Board, India).
\nocite{*}
\bibliographystyle{amsplain}
\bibliography{database}  
\end{document}